\documentclass[mn]{amsart}
\usepackage{amsfonts}

\usepackage{cite}

\newtheorem{theorem}{Theorem}
\theoremstyle{plain}

\newtheorem{corollary}[theorem]{Corollary}

\newtheorem{example}[theorem]{Example}

\newtheorem{proposition}[theorem]{Proposition}
\newtheorem{remark}[theorem]{Remark}

\numberwithin{equation}{section}

\input{tcilatex}

\begin{document}
\title[A Dirac delta operator]{{\huge A Dirac delta operator}}
\author{Juan Carlos Ferrando}
\address{Centro de Investigaci\'{o}n Operativa, Universidad 
\newline
\indent%
Miguel Hern\'{a}ndez, E-03202 Elche, Spain}
\email{jc.ferrando@umh.es}

\begin{abstract}
If $T$ is a (densely defined) self-adjoint operator acting on a complex
Hilbert space $\mathcal{H}$ and $I$ stands for the identity operator, we
introduce the delta function operator $\lambda \mapsto \delta \left( \lambda
I-T\right) $ at $T$. When $T$ is a bounded operator, then $\delta \left(
\lambda I-T\right) $ is an operator-valued distribution. If $T$ is
unbounded, $\delta \left( \lambda I-T\right) $ is a more general object that
still retains some properties of distributions. We derive various operative
formulas involving $\delta \left( \lambda I-T\right) $ and give several
applications of its usage.
\end{abstract}

\twocolumn[
\maketitle]

\section{The delta function $\protect\delta \left( \protect\lambda %
I-T\right) $}

The scalar delta `function' $\lambda \mapsto \delta \left( \lambda -a\right) 
$ along with its derivatives were introduced by Paul Dirac in \cite{Di1},
and later in \cite[Section 15]{Di2}, although its definition can be traced
back to Heaviside. The rigorous treatment of this object in the context of
distribution theory is due to Laurent Schwartz \cite{Ho, Sw3}. In this paper
we extend the definition of $\delta \left( \lambda -a\right) $ from real
numbers to self-adjoint operators on a Hilbert space $\mathcal{H}$. We
denote by $\mathcal{D}\left( \mathbb{R}\right) =\underrightarrow{\lim }\,%
\mathcal{D}\left( \left[ -n,n\right] \right) $ the linear space of
infinitely differentiable complex-valued functions of compact support,
equipped with the inductive limit topology. As usual in physics we shall
assume that the scalar product in $\mathcal{H}$ is anti-linear for the first
variable.

If $T$ is a densely defined self-adjoint operator\footnote{%
In what follows $\sigma \left( T\right) $ will denote the \textit{spectrum}
of $T$. Recall that the \textit{residual spectrum} of a self-adjoint
operator $T$ is empty, so that $\sigma \left( T\right) =\sigma _{p}\left(
T\right) \cup \sigma _{c}\left( T\right) $, where $\sigma _{p}\left(
T\right) $ denotes the \textit{point spectrum} (the eigenvalues) and $\sigma
_{c}\left( T\right) $ the \textit{continuous spectrum} of $T$.} on $\mathcal{%
H}$ and $I$ stands for the identity operator, we define the \textit{delta
function operator }$\lambda \mapsto \delta \left( \lambda I-T\right) $%
\textit{\ at }$T$ by%
\begin{equation}
f\left( T\right) =\int_{-\infty }^{+\infty }f\left( \lambda \right) \delta
\left( \lambda I-T\right) d\lambda  \label{Sym2}
\end{equation}%
for each $f\in \mathcal{C}\left( \mathbb{R}\right) $, i.\thinspace e., for
each real-valued continuous function $f\left( \lambda \right) $. Here $%
d\lambda $ is the Lebesgue measure of $\mathbb{R}$, but the right-hand side
of (\ref{Sym2}) is not a true integral. If $T$ is a bounded operator, we
shall see at once that $\delta \left( \lambda I-T\right) $ must be regarded
as a \textit{vector-valued distribution}, i.\thinspace e., as a continuous
linear map from the space $\mathcal{D}\left( \mathbb{R}\right) $ into the
locally convex space $\mathcal{L}\left( \mathcal{H}\right) $ of the bounded
linear operators (endomorphisms) on $\mathcal{H}$ equipped with the strong
operator topology \cite{Sw1, Sw2}, whose action on $f\in \mathcal{D}\left( 
\mathbb{R}\right) $ we denote as an integral. If $T$ is unbounded we shall
see that $\delta \left( \lambda I-T\right) $ still retains some useful
distributional-like properties. The previous equation means%
\begin{equation}
\left\langle y,f\left( T\right) x\right\rangle =\int_{-\infty }^{+\infty
}f\left( \lambda \right) \left\langle y,\delta \left( \lambda I-T\right)
x\right\rangle \,d\lambda  \label{weakeva}
\end{equation}%
for each $\left( x,y\right) \in D\left( f\left( T\right) \right) \times 
\mathcal{H}$, where $D\left( f\left( T\right) \right) $ stands for the
domain of the self-adjoint operator $f\left( T\right) $.

Let us recall that if $T$ is a (densely defined) self-adjoint operator,
there is a unique spectral family $\left\{ E_{\lambda }:\lambda \in \mathbb{R%
}\right\} $ of self-adjoint operators defined on the whole of $\mathcal{H}$
that satisfy $\left( i\right) $ $E_{\lambda }\leq E_{\mu }$ and $E_{\lambda
}E_{\mu }=E_{\lambda }$ for $\lambda \leq \mu $, $\left( ii\right) $ $%
\lim_{\epsilon \rightarrow 0^{+}}E_{\lambda +\epsilon }x=E_{\lambda }x$, and 
$\left( iii\right) $ $\lim_{\lambda \rightarrow -\infty }E_{\lambda }x=%
\mathbf{0}$ and $\lim_{\lambda \rightarrow \infty }E_{\lambda }x=x$ in $%
\mathcal{H}$ for all $x\in \mathcal{H}$. The domain $D\left( T\right) $ of $%
T $ consists of those $x\in \mathcal{H}$ such that%
\begin{equation*}
\int_{-\infty }^{+\infty }\left| \lambda \right| ^{2}\,d\left\| E_{\lambda
}x\right\| ^{2}<\infty .
\end{equation*}%
In this case, the spectral theorem (\textit{cf}. \cite[Section 107]{RN}) and
the Borel-measurable functional calculus provide a \textit{self-adjoint}
operator $f\left( T\right) $ defined by%
\begin{equation}
f\left( T\right) =\int_{-\infty }^{+\infty }f\left( \lambda \right)
\,dE_{\lambda }  \label{Sym1}
\end{equation}%
for each Borel-measurable function $f\left( \lambda \right) $, whose domain%
\begin{equation*}
D\left( f\left( T\right) \right) =\left\{ x\in \mathcal{H}:\int_{-\infty
}^{+\infty }\left| f\left( \lambda \right) \right| ^{2}\,d\left\| E_{\lambda
}x\right\| ^{2}<\infty \right\}
\end{equation*}%
is dense in $\mathcal{H}$. Observe that if $T$ is bounded, $f\left( T\right) 
$ need not be bounded. Moreover, since $\lambda \mapsto E_{\lambda }$ is
constant on the set $\mathbb{R}\setminus \sigma \left( T\right) $ of $T$, an
open set in $\mathbb{R}$, equation (\ref{Sym1}) tell us that $f\left(
\lambda \right) $ need not be defined on $\mathbb{R}\setminus \sigma \left(
T\right) $.

Thanks to (\ref{Sym1}) the definition of $\delta \left( \lambda I-T\right) $
may be extended to Borel-measurable functions by declaring that the equation
(\ref{Sym2}) holds for $\left( x,y\right) \in D\left( f\left( T\right)
\right) \times \mathcal{H}$ and each Borel function $f$. But, by reasons
that will become clear later, we shall restrict ourselves to those Borel
functions which are continuous at each point of $\sigma _{p}\left( T\right) $%
. Moreover, working with the real and complex parts, no difficulty arises if
the function $f$ involved in the equation (\ref{Sym2}) is complex-valued
(except that $f\left( T\right) $ is no longer a self-adjoint operator
whenever $\func{Im}f\neq 0$). Thus, unless otherwise stated, we shall assume
that both in (\ref{Sym2}) and (\ref{Sym1}) the function $f$ is
complex-valued. Note that the complex Stieltjes measure $d\left\langle
E_{\lambda }x,y\right\rangle $ need not be $d\lambda $-continuous. In what
follows we shall denote by $\mathcal{B}_{p}\left( \mathbb{R}\right) $ the
linear space over $\mathbb{C}$ consisting of all complex-valued
Borel-measurable functions of one real variable which are continuous on $%
\sigma _{p}\left( T\right) $.

If $f_{n}\rightarrow f$ in $\mathcal{D}\left( \mathbb{R}\right) $, the
sequence $\left\{ f_{n}\right\} _{n=1}^{\infty }$ is uniformly bounded and $%
f_{n}\left( x\right) \rightarrow f\left( x\right) $ at each $x\in \mathbb{R}$%
. So, if $T$ is bounded on $\mathcal{H}$ (equivalently, self-adjoint on the
whole of $\mathcal{H}$) it turns out that $f_{n}\left( T\right) \rightarrow
f\left( T\right) $ in the strong operator topology \cite[10.2.8 Corollary]%
{DS}. Therefore, in this case $\delta \left( \lambda I-T\right) $ is an $%
\mathcal{L}\left( \mathcal{H}\right) $-valued distribution.

As all integrals considered so far are over $\sigma \left( T\right) $, we
have 
\begin{equation}
\delta \left( \lambda I-T\right) =\mathbf{0}\ \ \forall \lambda \notin
\sigma \left( T\right) .  \label{Zero}
\end{equation}%
Also $\delta \left( -\lambda I+T\right) =\delta \left( \lambda I-T\right) $
for all $\lambda \in \mathbb{R}$. On the other hand, if $\mu \in \sigma
_{p}\left( T\right) $ and $y$ is an eigenvector corresponding to the
eigenvalue $\mu $, clearly%
\begin{equation}
\delta \left( \lambda I-T\right) y=\delta \left( \lambda -\mu \right) y
\label{Eigen}
\end{equation}%
for every $\lambda \in \mathbb{R}$. In the particular case when $T_{a}$ is
the linear operator defined on $\mathcal{H}$ by $T_{a}x=ax$ for a fixed $%
a\in \mathbb{R}$, then $T_{a}$ is a self-adjoint linear operator with $%
\sigma \left( T_{a}\right) =\sigma _{p}\left( T_{a}\right) =\left\{
a\right\} $. In this case $\delta \left( \lambda I-T_{a}\right) x=\delta
\left( \lambda -a\right) x$ for every $x\in \mathcal{H}$, i.\thinspace e., $%
\delta \left( \lambda I-T_{a}\right) =\delta \left( \lambda -a\right) I$.

Since equality $\left\langle f\left( T\right) ^{\dagger }y,x\right\rangle
=\left\langle y,f\left( T\right) x\right\rangle $ holds for all $x,y\in
D\left( f\left( T\right) \right) $ and each $f\in \mathcal{B}_{p}\left( 
\mathbb{R}\right) $, we may infer that%
\begin{equation*}
\left\langle y,\delta \left( \lambda I-T\right) x\right\rangle =\left\langle
\delta \left( \lambda I-T\right) y,x\right\rangle
\end{equation*}%
holds (in a `distributional' sense) for all $x,y\in D\left( T\right) $. This
suggests that in certain sense $\delta \left( \lambda I-T\right) $ may be
regarded (possibly for almost all $\lambda \in \mathbb{R}$) as a Hermitian
operator on $D\left( T\right) $.

Let us also point out that as equation (\ref{Sym2}) holds for all $f\in 
\mathcal{D}\left( \mathbb{R}\right) $, in a distributional sense we have%
\begin{equation}
\frac{d}{d\lambda }\left\langle y,E_{\lambda }x\right\rangle =\left\langle
y,\delta \left( \lambda I-T\right) x\right\rangle  \label{Sym3}
\end{equation}%
If $\lambda \mapsto Y\left( \lambda -\mu \right) $ denotes the unit step
function at $\mu \in \mathbb{R}$, given by $Y\left( \lambda -\mu \right) =0$
if $\lambda <\mu $ and $Y\left( \lambda -\mu \right) =1$ if $\lambda \geq
\mu $, since $E_{\lambda }=Y\left( \lambda I-T\right) $ for each $\lambda
\in \mathbb{R}$, formally%
\begin{equation}
dE_{\lambda }/d\lambda =Y^{\prime }\left( \lambda I-T\right) .  \label{Sym4}
\end{equation}%
So, from (\ref{Sym3}) and (\ref{Sym4}) we get $Y^{\prime }\left( \lambda
I-T\right) =\delta \left( \lambda I-T\right) $.

\begin{proposition}
\label{Differ}If $T$ is a bounded self-adjoint operator on $\mathcal{H}$ and 
$f\in \mathcal{C}^{1}\left( \mathbb{R}\right) $, then%
\begin{equation*}
\int_{-\infty }^{+\infty }f\left( \lambda \right) \,\delta ^{\prime }\left(
\lambda I-T\right) \,d\lambda =-f^{\prime }\left( T\right) .
\end{equation*}%
The same equality holds if $T$ is unbounded but $f\in \mathcal{D}\left( 
\mathbb{R}\right) $.
\end{proposition}

If $T$ is a self-adjoint operator and $f\in \mathcal{B}_{p}\left( \mathbb{R}%
\right) $, then%
\begin{equation*}
\int_{-\infty }^{+\,\infty }\left| f\left( \lambda \right) \right|
^{2}\delta \left( \lambda I-T\right) \,d\lambda =\int_{-\infty }^{+\,\infty
}\left| f\left( \lambda \right) \right| ^{2}\,dE_{\lambda }.
\end{equation*}%
where the latter equality \textit{is the definition} of $\left| f\left(
T\right) \right| ^{2}$. So, we have the following result.

\begin{proposition}
If $T$ is self-adjoint and $f\in \mathcal{B}_{p}\left( \mathbb{R}\right) $,
then%
\begin{equation*}
\left\langle f\left( T\right) y,f\left( T\right) x\right\rangle
=\int_{-\infty }^{+\,\infty }\left| f\left( \lambda \right) \right|
^{2}\left\langle y,\delta \left( \lambda I-T\right) x\right\rangle \,d\lambda
\end{equation*}%
for every $x,y\in D\left( f\left( T\right) \right) $.
\end{proposition}

\proof%
We adapt a classic argument. Indeed, for every $x,y\in D\left( f\left(
T\right) \right) $ we have%
\begin{equation*}
\left\langle f\left( T\right) y,f\left( T\right) x\right\rangle
=\int_{-\infty }^{+\,\infty }\overline{f\left( \lambda \right) }%
\,d\left\langle f\left( T\right) y,E_{\lambda }x\right\rangle .
\end{equation*}%
Since $E_{\mu }E_{\lambda }=E_{\mu }$ whenever $\mu \leq \lambda $, and $%
\left\langle E_{\lambda }y,x\right\rangle $ does not depend on $\mu $, by
splitting the integral we get 
\begin{equation}
\int_{-\infty }^{+\,\infty }f\left( \mu \right) \,d\left\langle E_{\mu
}y,E_{\lambda }x\right\rangle =\int_{-\infty }^{\lambda }f\left( \mu \right)
\,d\left\langle y,E_{\mu }x\right\rangle ,  \label{Torne}
\end{equation}%
where clearly the first integral is $\left\langle f\left( T\right)
y,E_{\lambda }x\right\rangle $. Plugging $d\left\langle f\left( T\right)
y,E_{\lambda }x\right\rangle $ into (\ref{Torne}), we are done. 
\endproof%

\begin{corollary}
Under the same conditions of the previous theorem, the equality%
\begin{equation}
\left\| f\left( T\right) x\right\| ^{2}=\int_{-\infty }^{+\,\infty }\left|
f\left( \lambda \right) \right| ^{2}\left\langle x,\delta \left( \lambda
I-T\right) x\right\rangle \,d\lambda  \label{Mary}
\end{equation}%
holds for every $x\in D\left( f\left( T\right) \right) $.
\end{corollary}

\begin{proposition}
\label{Angela}If $T$ is self-adjoint and $\left\{ f_{n}\right\}
_{n=1}^{\infty }$ is a uniformly bounded sequence in $\mathcal{B}_{p}\left( 
\mathbb{R}\right) $ such that $f_{n}\rightarrow f$ pointwise on $\mathbb{R}$
with $f\in \mathcal{B}_{p}\left( \mathbb{R}\right) $, then $f_{n}\left(
T\right) x\rightarrow f\left( T\right) x$ for every $x\in D\left( T\right) $.
\end{proposition}

\proof%
This is a straightforward consequence of preceding corollary and the
Lebesgue dominated convergence theorem. 
\endproof%

This proposition holds in particular if $f_{n}\rightarrow f$ in $\mathcal{D}%
\left( \mathbb{R}\right) $. Hence, even in the unbounded case, $\delta
\left( \lambda I-T\right) $ behaves as a vector-valued distribution-like
object.

\begin{proposition}
\label{Inter}Let $\left( \lambda ,\mu \right) \mapsto g\left( \lambda ,\mu
\right) $ be a function defined on $\mathbb{R}^{2}$ such that $g\left(
\lambda ,\cdot \right) \in \mathcal{L}_{1}\left( \mathbb{R}\right) $ for
every $\lambda \in \mathbb{R}$ and $g\left( \cdot ,\mu \right) \in \mathcal{B%
}_{p}\left( \mathbb{R}\right) $ for every $\mu \in \mathbb{R}$. If the
parametric integral 
\begin{equation*}
f\left( \lambda \right) =\int_{-\infty }^{+\,\infty }g\left( \lambda ,\mu
\right) d\mu
\end{equation*}%
is continuous on $\mathbb{R}$ and makes sense if we replace $\lambda $ by a
self-adjoint operator $T$, the value of the integral%
\begin{equation*}
\int_{-\infty }^{+\infty }\int_{-\infty }^{+\infty }g\left( \lambda ,\mu
\right) \delta \left( \lambda I-T\right) \,d\mu \,d\lambda
\end{equation*}%
does not depend on the integration ordering.
\end{proposition}

\proof%
Since $g\left( \cdot ,\mu \right) \in \mathcal{B}_{p}\left( \mathbb{R}%
\right) $ for every $\mu \in \mathbb{R}$., one has 
\begin{equation*}
g\left( T,\mu \right) =\int_{-\infty }^{+\infty }g\left( \lambda ,\mu
\right) \delta \left( \lambda I-T\right) d\lambda ,
\end{equation*}%
which implies%
\begin{equation*}
f\left( T\right) =\int_{-\infty }^{+\infty }\left\{ \int_{-\infty }^{+\infty
}g\left( \lambda ,\mu \right) \delta \left( \lambda I-T\right) d\lambda
\right\} d\mu .
\end{equation*}%
On the other hand, by the definition of $\delta \left( \lambda I-T\right) $
we have 
\begin{equation*}
f\left( T\right) =\int_{-\infty }^{+\infty }\left\{ \int_{-\infty }^{+\infty
}g\left( \lambda ,\mu \right) d\mu \right\} \delta \left( \lambda I-T\right)
\,d\lambda ,
\end{equation*}%
for $\left( x,y\right) \in D\left( T\right) \times \mathcal{H}$. So, the
proposition follows. 
\endproof%

\begin{theorem}
If $T$ is a self-adjoint operator on $\mathcal{H}$, then%
\begin{gather}
\int_{0}^{+\infty }f\left( \lambda \right) \delta \left( \lambda
I-T^{2}\right) d\lambda =  \label{Claire} \\
\int_{0}^{+\infty }\frac{1}{2\sqrt{\lambda }}\left\{ \delta \left( \sqrt{%
\lambda }I-T\right) -\delta \left( \sqrt{\lambda }I+T\right) \right\}
f\left( \lambda \right) d\lambda  \notag
\end{gather}%
if $\lambda >0$ and $f\in \mathcal{B}_{p}\left( \mathbb{R}\right) $, both
members acting on $D\left( T^{2}\right) $.
\end{theorem}

\proof%
First note that $T^{2}\geq 0$. Hence $\sigma \left( T^{2}\right) \subseteq %
\left[ 0,+\,\infty \right) $, which implies that $\delta \left( \lambda
I-T^{2}\right) =\mathbf{0}$ if $\lambda <0$. Since $T^{2}$ is a self-adjoint
operator, for $f\in \mathcal{B}_{p}\left( \mathbb{R}\right) $ we have%
\begin{equation*}
\int_{0}^{+\infty }f\left( \lambda \right) \delta \left( \lambda
I-T^{2}\right) d\lambda =f\left( T^{2}\right)
\end{equation*}%
On the other hand, it is clear that%
\begin{equation*}
\int_{0}^{+\infty }\frac{f\left( \lambda \right) }{2\sqrt{\lambda }}\,\delta
\left( \sqrt{\lambda }I-T\right) d\lambda =\int_{0}^{+\infty }f\left( \mu
^{2}\right) \delta \left( \mu I-T\right) d\mu
\end{equation*}%
whereas, using that $\delta \left( -\mu I+T\right) =\delta \left( \mu
I-T\right) $, we have%
\begin{equation*}
-\int_{0}^{+\infty }\frac{f\left( \lambda \right) }{2\sqrt{\lambda }}%
\,\delta \left( \sqrt{\lambda }I+T\right) d\lambda =\int_{-\infty
}^{0}f\left( \mu ^{2}\right) \delta \left( \mu I-T\right) d\mu
\end{equation*}%
So, the right-hand side of (\ref{Claire}) coincides with%
\begin{equation*}
\int_{-\infty }^{+\infty }f\left( \mu ^{2}\right) \delta \left( \mu
I-T\right) d\mu =f\left( T^{2}\right)
\end{equation*}%
since $\mu \mapsto f\left( \mu ^{2}\right) $ is a Borel function. 
\endproof%

If we denote by $L\left( \mathcal{H}\right) $ the linear space of all linear
endomorphisms on $\mathcal{H}$, the next theorem summarize some previous
results.

\begin{theorem}
If $T$ is a densely defined self-adjoint operator on a Hilbert space $%
\mathcal{H}$, there is an $L\left( \mathcal{H}\right) $-valued linear map $%
\delta _{T}$ on $\mathcal{B}_{p}\left( \mathbb{R}\right) $, whose action on $%
f\in \mathcal{B}_{p}\left( \mathbb{R}\right) $ we denote by%
\begin{equation*}
\left\langle \delta _{T},f\right\rangle =\int_{-\infty }^{+\,\infty }f\left(
\lambda \right) \,\delta \left( \lambda I-T\right) d\lambda ,
\end{equation*}%
such that $\left\langle \delta _{T},f\right\rangle =f\left( T\right) $. If $%
\left\{ f_{n}\right\} \subseteq \mathcal{B}_{p}\left( \mathbb{R}\right) $ is
uniformly bounded and $f_{n}\left( t\right) \rightarrow f\left( t\right) $,
with $f\in \mathcal{B}_{p}\left( \mathbb{R}\right) $, for all $t\in \mathbb{R%
}$ then $\left\langle \delta _{T},f_{n}\right\rangle x\rightarrow
\left\langle \delta _{T},f\right\rangle x$ for all $x\in \mathcal{H}$. If $T$
is bounded, $\delta _{T}$ is an $\mathcal{L}\left( \mathcal{H}\right) $%
-valued distribution, so $\left\langle \delta _{T},f\right\rangle $ is a
bounded operator on $\mathcal{H}$. In addition $\delta \left( \lambda
I-T\right) =\mathbf{0}$ if $\lambda \notin \sigma \left( T\right) $ and $%
\left\langle y,\delta \left( \lambda I-T\right) x\right\rangle =\left\langle
\delta \left( \lambda I-T\right) y,x\right\rangle $ for $x,y\in D\left(
T\right) $.
\end{theorem}

\section{Explicit form of $\protect\delta \left( \protect\lambda I-T\right) $%
}

If $Q$ is a vector-valued distribution, the \textit{Fourier transform} of $Q$
is defined as the vector valued distribution $\mathcal{F}Q$ on $\mathcal{S}%
\left( \mathbb{R}\right) $ such that $\left\langle \mathcal{F}%
Q,f\right\rangle =\left\langle Q,\mathcal{F}f\right\rangle $. As usual, we
denote by $\mathcal{F}^{-1}$ the inverse Fourier transform.

\begin{theorem}
\label{Furie}If $T$ is a self-adjoint operator, the identity%
\begin{equation}
\delta \left( \lambda I-T\right) =\frac{1}{2\pi }\int_{-\infty }^{+\,\infty
}e^{it\left( \lambda I-T\right) }dt  \label{Expo}
\end{equation}%
holds for every $\lambda \in \mathbb{R}$, and the action $f\left( T\right) $
of $\delta \left( \lambda I-T\right) $ on $f\in \mathcal{S}\left( \mathbb{R}%
\right) $ is given by%
\begin{equation*}
f\left( T\right) =\int_{-\infty }^{+\infty }\left\{ \int_{-\infty }^{+\infty
}\frac{f\left( \lambda \right) }{2\pi }\,e^{it\left( \lambda I-T\right)
}\,d\lambda \right\} dt.
\end{equation*}
\end{theorem}

\proof%
Setting $\delta _{T}\left( \lambda \right) =\delta \left( \lambda I-T\right) 
$ observe that%
\begin{equation*}
\left( \mathcal{F}\delta _{T}\right) \left( t\right) =\frac{1}{\sqrt{2\pi }}%
\,e^{-itT}.
\end{equation*}%
Indeed, if $f\in \mathcal{S}\left( \mathbb{R}\right) $ we have%
\begin{eqnarray*}
\left\langle \mathcal{F}\delta _{T},f\right\rangle &=&\left\langle \delta
_{T},\mathcal{F}f\right\rangle =\int_{-\infty }^{+\,\,\infty }\left( 
\mathcal{F}f\right) \left( \lambda \right) \delta \left( \lambda I-T\right)
d\lambda \\
&=&\left( \mathcal{F}f\right) \left( T\right) =\frac{1}{\sqrt{2\pi }}%
\int_{-\infty }^{+\,\infty }f\left( t\right) e^{-itT}dt.
\end{eqnarray*}%
Consequently%
\begin{equation}
\delta _{T}=\mathcal{F}^{-1}\left\{ \frac{1}{\sqrt{2\pi }}\,e^{-itT}\right\}
.  \label{Minerva}
\end{equation}%
Functionally, the action of $\delta _{T}$\ on $f\in S\left( \mathbb{R}%
\right) $\ by means of equation (\ref{Minerva}) becomes%
\begin{equation}
\left\langle \delta _{T},f\right\rangle =\left\langle \frac{1}{\sqrt{2\pi }}%
\,e^{-itT},\left( \mathcal{F}^{-1}f\right) \left( t\right) \right\rangle
\label{Fortune}
\end{equation}%
Consequently, we have%
\begin{equation*}
\left\langle \delta _{T},f\right\rangle =\int_{-\infty }^{+\infty }\left\{
\int_{-\infty }^{+\infty }\frac{f\left( \mu \right) }{2\pi }\,e^{it\left(
\mu I-T\right) }\,d\mu \right\} dt
\end{equation*}%
with the order of the integration as stated. 
\endproof%

\begin{corollary}
\label{Order}If $e^{-itT}x=x\left( t\right) $, for $x\in D\left( T\right) $
one has%
\begin{equation*}
\delta \left( \lambda I-T\right) x=\frac{1}{2\pi }\int_{-\infty }^{+\infty
}e^{i\lambda t}x\left( t\right) \,dt
\end{equation*}%
and if $x\in D\left( f\left( T\right) \right) $ and $f\in \mathcal{S}\left( 
\mathbb{R}\right) $, then%
\begin{equation*}
f\left( T\right) x=\frac{1}{2\pi }\int_{-\infty }^{+\infty }\left\{
\int_{-\infty }^{+\infty }f\left( \lambda \right) e^{i\lambda t}d\lambda
\right\} x\left( t\right) \,dt.
\end{equation*}
\end{corollary}

\begin{remark}
\emph{Consider the one-parameter unitary group }$\left\{ U\left( t\right)
:t\in \mathbb{R}\right\} $\emph{\ generated by the self-adjoint operator }$T$%
\emph{, that is, }$U\left( t\right) =\exp \left( -itT\right) $\emph{\ for
every }$t\in \mathbb{R}$\emph{. If} $\mathcal{F}$ \emph{denotes the Fourier
transform, equation (\ref{Minerva}) can be written as}%
\begin{equation}
\delta \left( \lambda I-T\right) =\frac{1}{\sqrt{2\pi }}\,\mathcal{F}%
^{-1}\left( U\right) \left( \lambda \right) .  \label{Zeus1}
\end{equation}%
\emph{So, equation (\ref{Fortune}) reads as}%
\begin{equation}
f\left( T\right) =\frac{1}{\sqrt{2\pi }}\int_{-\infty }^{+\infty }\left( 
\mathcal{F}^{-1}f\right) \left( t\right) U\left( t\right) \,dt.
\label{Zeus2}
\end{equation}
\end{remark}

In what follows we shall compute the spectral family $\left\{ E_{\lambda
}:\lambda \in \mathbb{R}\right\} $ for some useful self-adjoint operators of
Quantum Mechanics by means of the delta $\delta \left( \lambda I-T\right) $.
Nonetheless, although $E_{\lambda }=Y\left( \lambda I-T\right) $, the
identification 
\begin{equation*}
Y\left( \lambda I-T\right) =\int_{-\infty }^{+\,\infty }Y\left( \lambda -\mu
\right) \delta \left( \mu I-T\right) d\mu
\end{equation*}%
might be not well-defined because $\mu \mapsto Y\left( \lambda -\mu \right) $
has a jump discontinuity at $\mu =\lambda $. Indeed, if $\lambda \in \sigma
_{p}\left( T\right) $ and $x$ is an eigenvector corresponding to $\lambda $,
then%
\begin{equation*}
\int_{-\infty }^{+\,\infty }Y\left( \lambda -\mu \right) \delta \left( \mu
I-T\right) x\,d\mu =\left\{ \int_{-\infty }^{\lambda }\delta \left( \mu
-\lambda \right) \,d\mu \right\} x
\end{equation*}%
and the right-hand integral makes no sense (see \cite{GW} for a useful
discussion). If $\lambda \notin \sigma _{p}\left( T\right) $ we define%
\begin{equation}
E_{\lambda }=\int_{-\infty }^{+\,\infty }Y\left( \lambda -\mu \right) \delta
\left( \mu I-T\right) d\mu  \label{Hamlet}
\end{equation}%
If $\lambda $ belongs to $\sigma _{p}\left( T\right) $, then $\left( \mu
\mapsto Y\left( \lambda -\mu \right) \right) \notin \mathcal{B}_{p}\left( 
\mathbb{R}\right) $. In order to define $E_{\lambda }$ we enlarge a little
the interval of integration by considering the integral%
\begin{equation*}
\int_{-\infty }^{\lambda +\epsilon }\delta \left( \mu -\lambda \right) \,d\mu
\end{equation*}%
for small $\epsilon >0$. So, if $\lambda \in \sigma _{p}\left( T\right) $ we
define%
\begin{equation}
E_{\lambda }=\lim_{\epsilon \rightarrow 0^{+}}\int_{-\infty }^{+\,\infty
}Y\left( \lambda +\epsilon -\mu \right) \delta \left( \mu I-T\right) d\mu .
\label{Macbeth}
\end{equation}%
The limit is well-defined since $\lim_{\epsilon \rightarrow 0^{+}}E_{\lambda
+\epsilon }=E_{\lambda }$ pointwise on $\mathcal{H}$. In the particular case
when $\lambda $ belongs to $\sigma _{d}\left( T\right) $, the discrete part
of $\sigma _{p}\left( T\right) $, $\lambda $ is isolated in $\sigma
_{p}\left( T\right) $.

\begin{example}
\label{Hugh}The spectral family of the \emph{(}up to a sign\emph{)}
one-dimensional Quantum Mechanics momentum operator of the free particle $%
P=iD$, where $D\varphi =\varphi ^{\prime }$, acting on the Hilbert space $%
\mathcal{H}=L_{2}\left( \mathbb{R}\right) $ is given by%
\begin{equation*}
\left( E_{\lambda }\varphi \right) \left( x\right) =\frac{1}{2}\,\varphi
\left( x\right) +\frac{1}{2\pi i}\,\mathrm{p.v.}\int_{-\infty }^{+\infty }%
\frac{e^{i\lambda \left( s-x\right) }}{s-x}\,\varphi \left( s\right) ds
\end{equation*}%
for every regular compactly supported $\varphi \in D\left( P\right) $.
\end{example}

\proof%
As is well-known $P$ is a self-adjoint operator with $D\left( P\right)
=H^{2,1}\left( \mathbb{R}\right) $ and $\sigma _{c}\left( P\right) =\mathbb{R%
}$. Since%
\begin{equation*}
\left( e^{-itP}\varphi \right) \left( x\right) =\left( e^{tD}\varphi \right)
\left( x\right) =\varphi \left( x+t\right)
\end{equation*}%
for a regular enough $\varphi \in D\left( P\right) $, by Corollary \ref%
{Order} we have%
\begin{equation*}
\left\{ \delta \left( \mu I-P\right) \varphi \right\} \left( x\right) =\frac{%
1}{2\pi }\int_{-\infty }^{+\infty }e^{i\mu t}\varphi \left( x+t\right) dt.
\end{equation*}%
Note that the integral of the right-hand side does exist because $\varphi $
has compact support.

According to the definition of $E_{\lambda }$ for the continuous spectrum
and keeping in mind the order of integration as indicated in Corollary \ref%
{Order}, one has%
\begin{equation*}
\left\{ E_{\lambda }\varphi \right\} \left( x\right) =\frac{1}{2\pi }%
\int_{-\infty }^{+\,\infty }\int_{-\infty }^{+\,\infty }Y\left( \lambda -\mu
\right) e^{i\mu t}\varphi \left( x+t\right) \,d\mu \,dt.
\end{equation*}%
So, since%
\begin{equation*}
\frac{1}{\sqrt{2\pi }}\int_{-\infty }^{+\,\infty }Y\left( \lambda -\mu
\right) e^{i\mu t}d\mu =\mathcal{F}\left( Y\right) \left( t\right) \cdot
e^{i\lambda t},
\end{equation*}%
bearing in mind the distributional relation%
\begin{equation}
\mathcal{F}\left( Y\right) \left( t\right) =\sqrt{\frac{\pi }{2}}\left(
\delta \left( t\right) +\frac{1}{i\pi }\,\mathrm{p.v.}\frac{1}{t}\right) ,
\label{Menge}
\end{equation}%
we get 
\begin{equation*}
\left\{ E_{\lambda }\varphi \right\} \left( x\right) =\frac{1}{2}\,\varphi
\left( x\right) +\frac{1}{2\pi i}\int_{-\infty }^{+\infty }\frac{e^{i\lambda
\left( s-x\right) }}{s-x}\,\varphi \left( s\right) ds
\end{equation*}%
where the last integral must be understood in Cauchy's principal value
sense. 
\endproof%

\begin{example}
The spectral family of the one-dimensional Quantum Mechanics kinetic energy
term of the free particle, corresponding to the Laplace operator $T=-D^{2}$
on $\mathcal{H}=L_{2}\left( \mathbb{R}\right) $, where $D^{2}\varphi
=\varphi ^{\prime \prime }$, is given by%
\begin{equation*}
\left( E_{\lambda }\varphi \right) \left( x\right) =\frac{1}{i\pi }\,\mathrm{%
p.v.}\int_{-\infty }^{+\infty }\frac{\cos \left( \lambda \left( s-x\right)
\right) -1}{s-x}\,\varphi \left( s\right) ds
\end{equation*}%
for $\lambda >0$ and $E_{\lambda }=\mathbf{0}$ whenever $\lambda <0$, where $%
\varphi $ is a regular function with compact support belonging to $D\left(
T\right) $.
\end{example}

\proof%
In this case $T$ is a self-adjoint operator with $\sigma \left( T\right) =%
\left[ 0,+\infty \right) $. Since $T=\left( iD\right) ^{2}$, according to (%
\ref{Claire}) we have%
\begin{equation*}
\delta \left( \lambda I-T\right) =\frac{1}{2\sqrt{\lambda }}\left\{ \delta
\left( \sqrt{\lambda }I-iD\right) -\delta \left( \sqrt{\lambda }I+iD\right)
\right\}
\end{equation*}%
regarded as a functional on $\mathcal{S}\left( \mathbb{R}\right) $ through $%
d\lambda $-integration over $\left[ 0,+\infty \right) $. Plugging%
\begin{equation*}
\left( \delta \left( \mu I\mp iD\right) \varphi \right) \left( x\right) =%
\frac{1}{2\pi }\int_{-\infty }^{+\infty }e^{i\mu t}\varphi \left( x\pm
t\right) dt
\end{equation*}%
into the previous expression and keeping in mind the correct order of
integration, we see that%
\begin{gather*}
\int_{0}^{\infty }f\left( \lambda \right) \left( \delta \left( \lambda
I-T\right) \varphi \right) \left( x\right) d\lambda = \\
\frac{1}{4\pi }\int_{0}^{\infty }\int_{-\infty }^{+\infty }f\left( \lambda
\right) \frac{e^{i\sqrt{\lambda }t}}{\sqrt{\lambda }}\left[ \,\varphi \left(
x+t\right) -\varphi \left( x-t\right) \right] \,d\lambda \,dt
\end{gather*}%
for every $f\in \mathcal{S}\left( \mathbb{R}\right) $. By the definition of $%
E_{\lambda }$ if $\lambda >0$ and the fact that $\delta \left( \mu
I-T\right) =\mathbf{0}$ whenever $\mu <0$, we have%
\begin{equation*}
\left( E_{\lambda }\varphi \right) \left( x\right) =\int_{0}^{+\,\infty
}Y\left( \lambda -\mu \right) \left( \delta \left( \mu I-T\right) \varphi
\right) \left( x\right) d\mu .
\end{equation*}%
Working out the penultimate integral with $\mu $ instead of $\lambda $ and $%
f\left( \mu \right) =Y\left( \lambda -\mu \right) $, we obtain 
\begin{gather*}
\int_{0}^{+\,\infty }\int_{-\infty }^{+\infty }Y\left( \lambda -\mu \right) 
\frac{e^{i\sqrt{\mu }t}}{\sqrt{\mu }}\left[ \,\varphi \left( x+t\right)
-\varphi \left( x-t\right) \right] \,d\mu \,dt \\
=\int_{-\infty }^{+\infty }\left\{ \int_{0}^{\lambda }\frac{e^{i\sqrt{\mu }t}%
}{\sqrt{\mu }}\,d\mu \right\} \left[ \,\varphi \left( x+t\right) -\varphi
\left( x-t\right) \right] \,dt
\end{gather*}%
for $\lambda >0$. So, by setting $u=\sqrt{\mu }$ we get%
\begin{equation*}
\left( E_{\lambda }\varphi \right) \left( x\right) =\int_{-\infty }^{+\infty
}\frac{dt}{2\pi }\left[ \,\varphi \left( x+t\right) -\varphi \left(
x-t\right) \right] \int_{0}^{\lambda }e^{iut}du.
\end{equation*}%
Now we have 
\begin{equation*}
\frac{1}{\sqrt{2\pi }}\int_{0}^{\lambda }e^{iut}du=\left( 1-e^{i\lambda
t}\right) \mathcal{F}^{-1}\left( Y\right) \left( t\right) ,
\end{equation*}%
so, using that $\mathcal{F}^{-1}\left( Y\left( v\right) \right) =\mathcal{F}%
\left( 1-Y\left( v\right) \right) \left( t\right) $ as well as equation (\ref%
{Menge}), we get%
\begin{equation*}
\frac{1}{\sqrt{2\pi }}\int_{0}^{\lambda }e^{iut}du=\left( 1-e^{i\lambda
t}\right) \sqrt{\frac{\pi }{2}}\left( \delta \left( t\right) -\frac{1}{i\pi }%
\,\mathrm{p.v.}\frac{1}{t}\right)
\end{equation*}%
which implies%
\begin{gather*}
\left( E_{\lambda }\varphi \right) \left( x\right) = \\
-\frac{1}{\pi i}\int_{-\infty }^{+\infty }\frac{\varphi \left( s\right) }{s-x%
}\,ds+\frac{1}{\pi i}\int_{-\infty }^{+\infty }\frac{\cos \left( \lambda
\left( s-x\right) \right) }{s-x}\,\varphi \left( s\right) ds
\end{gather*}%
where the integrals are understood in Cauchy's principal value sense. 
\endproof%

\begin{example}
Spectral family of the \emph{(}up to a sign\emph{)} one-dimensional Quantum
Mechanics momentum operator $S$ for a bounded particle on $\mathcal{H}=L_{2}%
\left[ -\pi ,\pi \right] $ with domain%
\begin{equation*}
\left\{ \varphi \in L_{2}\left[ -\pi ,\pi \right] :\varphi ^{\prime }\in
L_{2}\left[ -\pi ,\pi \right] ,\,\varphi \left( -\pi \right) =\varphi \left(
\pi \right) \right\}
\end{equation*}%
\emph{As is well-known this is a self-adjoint operator with discrete
spectrum }$\sigma \left( S\right) =\mathbb{Z}$ \emph{whose eigenfunction
system }$\{\varphi _{n}:n\in \mathbb{Z}\}$\emph{, with }$\varphi _{n}\left(
x\right) =\left( 2\pi \right) ^{-1/2}e^{-inx}$\emph{, are the solutions of
the eigenvalue problem }$i\varphi ^{\prime }=\lambda \varphi $\emph{\ with }$%
\varphi \left( -\pi \right) =\varphi \left( \pi \right) $\emph{. So, for }$%
\varphi \in D\left( S\right) $\emph{\ we have }$\varphi \overset{L_{2}}{=}%
\tsum_{n\in \mathbb{Z}}c_{n}\varphi _{n}$ \emph{with}%
\begin{equation*}
c_{n}=\left\langle \varphi ,\varphi _{n}\right\rangle =\frac{1}{2\pi }%
\int_{-\pi }^{\pi }\varphi \left( x\right) e^{inx}dx
\end{equation*}%
\emph{for every }$n\in \mathbb{Z}$\emph{. Since }$\sigma \left( S\right)
=\sigma _{d}\left( S\right) $\emph{, recalling the definition of the
operator }$E_{\lambda }$\emph{\ for }$\lambda \in \sigma _{d}\left( S\right) 
$\emph{, clearly we have}%
\begin{equation*}
\left( E_{\lambda }\varphi \right) \left( x\right) =\lim_{\epsilon
\rightarrow 0^{+}}\int_{-\infty }^{+\,\infty }Y\left( \lambda +\epsilon -\mu
\right) \left( \delta \left( \mu I-S\right) \varphi \right) \left( x\right)
d\mu
\end{equation*}%
\emph{for every }$\lambda \in \mathbb{R}$\emph{. So, the fact that }$%
E_{\lambda }$\emph{\ is a bounded operator yields}%
\begin{equation*}
E_{\lambda }\varphi =\sum_{n\in \mathbb{Z}}c_{n}\,E_{\lambda }\varphi _{n}
\end{equation*}%
\emph{Using that }$\delta \left( \mu I-S\right) e^{-inx}=\delta \left( \mu
-n\right) e^{-inx}$\emph{\ and that }$Y\left( \lambda +0-n\right) =Y\left(
\lambda -n\right) $\emph{, we get} 
\begin{equation*}
\left( E_{\lambda }\varphi \right) \left( x\right) =\sum_{n\in \mathbb{Z}}%
\frac{c_{n}}{\sqrt{2\pi }}\,Y\left( \lambda -n\right) e^{-inx}=\sum_{n\in 
\mathbb{Z},\,n\leq \left[ \lambda \right] }\frac{e^{-inx}}{\sqrt{2\pi }}.
\end{equation*}
\end{example}

\begin{remark}
\emph{Since in the previous example }$S$\emph{\ is bounded on }$\mathcal{H}%
=L_{2}\left[ -\pi ,\pi \right] $\emph{, the delta operator }$\delta \left(
\lambda I-S\right) $\emph{\ should be regarded as a continuous endomorphism
as well. In this case }%
\begin{equation*}
\delta \left( \lambda I-S\right) \varphi =\sum_{n\in \mathbb{Z}%
}c_{n}\,\delta \left( \lambda -n\right) \varphi _{n}.
\end{equation*}
\end{remark}

\begin{example}
The one-dimensional Quantum Mechanics position operator on $L_{2}\left( 
\mathbb{R}\right) $. \emph{This operator is defined on }$\mathcal{H}%
=L_{2}\left( \mathbb{R}\right) $\emph{\ by }$\left( Q\varphi \right) \left(
x\right) =x\varphi \left( x\right) $\emph{\ for every }$x\in \mathbb{R}$%
\emph{. Clearly }$\sigma _{c}\left( Q\right) =\mathbb{R}$ \emph{and }$%
\varphi \in D\left( Q\right) $\emph{\ if }$\left( x\mapsto x\,\varphi \left(
x\right) \right) \in \mathcal{L}_{2}\left( \mathbb{R}\right) $\emph{.
Moreover, it is clear that}%
\begin{equation*}
\left\{ \exp \left( it\left( \lambda I-Q\right) \right) \varphi \right\}
\left( x\right) =e^{i\left( \lambda -x\right) t}\varphi \left( x\right) .
\end{equation*}%
\emph{So we have }%
\begin{equation*}
\left( \delta \left( \lambda I-Q\right) \varphi \right) \left( x\right)
=\delta \left( \lambda -x\right) \varphi \left( x\right) .
\end{equation*}%
\emph{Hence, in this case we can write}%
\begin{equation*}
\left\{ E_{\lambda }\varphi \right\} \left( x\right) =\int_{-\infty
}^{+\infty }Y\left( \lambda -\mu \right) \delta \left( \mu -x\right) \varphi
\left( x\right) d\mu
\end{equation*}%
\emph{Therefore, if }$\lambda \neq x$\emph{\ we get}%
\begin{equation*}
\left\{ E_{\lambda }\varphi \right\} \left( x\right) =Y\left( \lambda
-x\right) \varphi \left( x\right) .
\end{equation*}
\end{example}

\begin{example}
\label{Viola}Explicit form of $\delta \left( \lambda I-M\right) $ for the
Hermitian matrix of $\mathcal{H}=\mathbb{C}^{3}$%
\begin{equation*}
M=\left[ 
\begin{array}{rrr}
0 & 1 & 1 \\ 
1 & 0 & 1 \\ 
1 & 1 & 0%
\end{array}%
\right] .
\end{equation*}
\end{example}

\proof%
In this case $M=PJ_{M}P^{-1}$ with $\sigma \left( M\right) =\left\{
-1,2\right\} $ and%
\begin{equation*}
J_{M}=\left[ 
\begin{array}{rrr}
-1 & 0 & 0 \\ 
0 & -1 & 0 \\ 
0 & 0 & 2%
\end{array}%
\right] ,\quad P=\left[ 
\begin{array}{rrr}
1 & 1 & 1 \\ 
-1 & 0 & 1 \\ 
0 & -1 & 1%
\end{array}%
\right]
\end{equation*}%
Using (\ref{Expo}) together with the fact that%
\begin{equation*}
\dint_{-\infty }^{+\infty }e^{i\left( \lambda -\rho \right) t}dt=2\pi \delta
\left( \lambda -\rho \right) ,
\end{equation*}%
we get%
\begin{gather*}
\delta \left( \lambda I-M\right) =\frac{1}{2\pi }\,P\left\{ \int_{-\infty
}^{+\infty }\exp it\left( \lambda I-J_{M}\right) dt\right\} P^{-1} \\
=P\left[ 
\begin{array}{ccc}
\delta \left( \lambda +1\right) & 0 & 0 \\ 
0 & \delta \left( \lambda +1\right) & 0 \\ 
0 & 0 & \delta \left( \lambda -2\right)%
\end{array}%
\right] P^{-1}
\end{gather*}%
Let us compute the spectral family and the projection operator onto the
eigenspace $\ker \left( M+I\right) $. Clearly 
\begin{equation*}
E_{\lambda }=P\left[ 
\begin{array}{ccc}
Y\left( \lambda +1\right) & 0 & 0 \\ 
0 & Y\left( \lambda +1\right) & 0 \\ 
0 & 0 & Y\left( \lambda -2\right)%
\end{array}%
\right] P^{-1}
\end{equation*}%
for every $\lambda \in \mathbb{R}$. If $\lambda _{1}=-1$, the orthogonal
projection $P_{\lambda _{1}}$ onto $\ker \left( I+M\right) $ is%
\begin{equation*}
P_{\lambda _{1}}=\frac{1}{3}\,P\left[ 
\begin{array}{ccc}
1 & 0 & 0 \\ 
0 & 1 & 0 \\ 
0 & 0 & 0%
\end{array}%
\right] P^{-1}=\frac{1}{3}\left[ 
\begin{array}{rrr}
2 & -1 & -1 \\ 
-1 & 2 & -1 \\ 
-1 & -1 & 2%
\end{array}%
\right]
\end{equation*}%
since $P_{\lambda _{1}}=E_{\lambda _{1}}-E_{\lambda _{1}-0}=E_{\lambda _{1}}$%
. 
\endproof%

\begin{example}
Consider a compact self-adjoint operator $K$ acting on a separable Hilbert
space $\mathcal{H}$ which does not admit the eigenvalue zero. Let $\left\{
u_{i}:i\in \mathbb{N}\right\} $ be a Hilbert basis of $\mathcal{H}$ with its
corresponding sequence of real eigenvalues $\left\{ \lambda _{i}:i\in 
\mathbb{N}\right\} $, where $\left| \lambda _{i+1}\right| \leq \left|
\lambda _{i}\right| $ for every $i\in \mathbb{N}$. Let us compute the action
of the operator $\left( \lambda I-K\right) ^{-1}$ on any $x\in \mathcal{H}$
and the operator $\delta \left( \lambda I-T\right) $.
\end{example}

\proof%
If $x\in \mathcal{H}$, we can write $x=\sum_{i=1}^{\infty }\left\langle
x,u_{i}\right\rangle u_{i}$. Since $\left( \lambda I-K\right) ^{-1}$ is a
bounded operator whenever $\lambda \notin \sigma \left( K\right) $, we have%
\begin{equation*}
\left( \lambda I-K\right) ^{-1}x=\sum_{i=1}^{\infty }\left\langle
x,u_{i}\right\rangle \int_{-\infty }^{+\,\infty }\frac{1}{\lambda -\mu }%
\,\delta \left( \mu I-K\right) u_{i}
\end{equation*}%
so we obtain the classic series%
\begin{equation*}
\left( \lambda I-K\right) ^{-1}x=\sum_{i=1}^{\infty }\frac{1}{\lambda -\mu
_{i}}\,\left\langle x,u_{i}\right\rangle u_{i}.
\end{equation*}%
For the solution of the equation $\left( I-zK\right) x=y$\ with $z\in 
\mathbb{C}$ we get the Schmidt series%
\begin{equation*}
x=\left( I-zK\right) ^{-1}y=\sum_{i=1}^{\infty }\frac{1}{1-z\mu _{i}}%
\,\left\langle y,u_{i}\right\rangle u_{i}
\end{equation*}%
whenever $z^{-1}\notin \sigma \left( T\right) $. On the other hand, since $%
\delta \left( \lambda I-K\right) $ acts on $\mathcal{H}$ as a continuous
endomorphism, equation 
\begin{equation*}
\delta \left( \lambda I-K\right) x=\tsum_{i=1}^{\infty }\left\langle
x,u_{i}\right\rangle \delta \left( \lambda -\mu _{i}\right) u_{i}.
\end{equation*}%
holds for every $x\in \mathcal{H}$. 
\endproof%

If $T$ is an unbounded self-adjoint operator then $D\left( T\right) \neq 
\mathcal{H}$ and $D\left( T^{n}\right) $ becomes smaller as $n$ grows. So,
the following result, makes sense only if the operator $T$ is bounded.

\begin{theorem}
In general, if $T$ is a bounded self-adjoint operator, one has%
\begin{equation}
\delta \left( \lambda I-T\right) =\sum_{n=0}^{\infty }\left( -1\right) ^{n}%
\frac{\delta ^{\left( n\right) }\left( \lambda \right) }{n!}\,T^{n}
\label{Blondie}
\end{equation}%
which is the Taylor series of $\delta \left( \lambda I-T\right) $\ at $%
\lambda I$.
\end{theorem}

\proof%
Developing the\ operator function $\exp \left( itT\right) $, which is
well-defined by the spectral theorem, we get%
\begin{equation*}
\delta \left( \lambda I-T\right) =\frac{1}{2\pi }\int_{-\infty }^{+\infty
}e^{-i\lambda t}\sum_{n=0}^{\infty }\frac{\left( it\right) ^{n}}{n!}%
\,T^{n}dt,
\end{equation*}%
so that, formally interchanging the sum and the integral, we may write%
\begin{equation*}
\delta \left( \lambda I-T\right) =\frac{1}{\sqrt{2\pi }}\sum_{n=0}^{\infty }%
\frac{\mathcal{F}\left\{ \left( it\right) ^{n}\right\} \left( \lambda
\right) }{n!}\,T^{n}.
\end{equation*}%
Using the fact that%
\begin{equation*}
\mathcal{F}\left\{ \left( it\right) ^{n}\right\} \left( \lambda \right)
=\left( -1\right) ^{n}\sqrt{2\pi }\,\delta ^{\left( n\right) }\left( \lambda
\right)
\end{equation*}%
for every $n\in \mathbb{N}$, we obtain (\ref{Blondie}). 
\endproof%

\section{The resolvent operator and $\protect\delta \left( \protect\lambda %
I-T\right) $}

Recall that the spectrum $\sigma \left( T\right) $ of a (densely defined)
self-adjoint operator on a complex Hilbert space $\mathcal{H}$ is a closed
subset of $\mathbb{C}$ contained in $\mathbb{R}$ (see for instance \cite[3.2]%
{Sc}). If $z\in \mathbb{C}\setminus \sigma \left( T\right) $, i.\thinspace
e., if $z$ is a regular point of $T$, and%
\begin{equation*}
\mathcal{R}\left( z,T\right) =\left( zI-T\right) ^{-1}
\end{equation*}%
denotes the \textit{resolvent} operator of $T$ at $z$ (see \cite[Definition
8.2]{Mo}), the function $\lambda \mapsto \left( z-\lambda \right) ^{-1}$ is
continuous on $\sigma \left( T\right) $. The resolvent is well-defined over $%
\mathcal{H}$, so it is a bounded normal operator. If $z\in \mathbb{R}%
\setminus \sigma \left( T\right) $ then $\mathcal{R}\left( z,T\right) $ is
even self-adjoint. From (\ref{Sym2}) it follows that%
\begin{equation*}
\mathcal{R}\left( z,T\right) =\int_{-\infty }^{+\infty }\frac{1}{z-\lambda }%
\,\delta \left( \lambda I-T\right) \,d\lambda
\end{equation*}%
which is the \textit{integral form of the resolvent} of $T$. So, by
considering the complex-valued function $f\left( \lambda \right) =\left(
z-\lambda \right) ^{-1}$ with $z\in \mathbb{C}\setminus \sigma \left(
T\right) $ and using the fact that%
\begin{equation*}
\mathcal{F}^{-1}\left( \frac{1}{\lambda -z}\right) \left( t\right) =\sqrt{%
2\pi }\,ie^{izt}Y\left( t\right)
\end{equation*}%
then, according to (\ref{Zeus2}), for $\func{Im}z>0$ we have%
\begin{equation*}
\left( zI-T\right) ^{-1}=-i\int_{0}^{\infty }e^{izt}U\left( t\right) .
\end{equation*}%
From here, it follows that%
\begin{equation*}
\mathcal{R}\left( z,iT\right) =i\mathcal{R}\left( iz,-T\right) =\left( 
\mathcal{L}U^{-1}\right) \left( z\right)
\end{equation*}%
if $\func{Im}z>0$, where $\mathcal{L}$ is the Laplace transform. This is 
\textit{the Hille-Yosida theorem} which relates the resolvent with the
one-parameter group of unitary transformations $\left\{ U\left( t\right)
:t\in \mathbb{R}\right\} $ generated by the self-adjoint operator $T$.

If $T$ is a bounded self-adjoint operator, $\gamma $ is a closed Jordan
contour that encloses $\sigma \left( T\right) $ and $f\left( z\right) $ is
holomorphic inside the connected region surrounded by the path $\gamma $,
the Dunford integral formula asserts that%
\begin{equation*}
\frac{1}{2\pi i}\int_{\gamma }f\left( z\right) \mathcal{R}\left( z,T\right)
\,dz=f\left( T\right) .
\end{equation*}%
In \cite{Ti} is pointed out that $\left( 2\pi i\right) ^{-1}\mathcal{R}%
\left( z,T\right) $ can be considered as the indicatrix of a vector-valued
distribution with values in $\mathcal{L}\left( \mathcal{H}\right) $. Dunford
integral formula is easily obtained by using the $\delta \left( \lambda
I-T\right) $ operator since, if we apply the Proposition \ref{Inter} with $%
g\left( \lambda ,\mu \right) =f\left( z\left( \mu \right) \right) \left(
z\left( \mu \right) -\lambda \right) ^{-1}$, where $z\left( \mu \right)
=\gamma \left( \mu \right) $ and $0\leq \mu \leq 1$, then%
\begin{gather*}
\int_{\gamma }f\left( z\right) \mathcal{R}\left( z,T\right)
dz\,=\int_{-\infty }^{+\infty }\left\{ \int_{\gamma }\frac{f\left( z\right) 
}{z-\lambda }\,dz\right\} \delta \left( \lambda I-T\right) d\lambda \\
=2\pi i\int_{-\infty }^{+\infty }f\left( \lambda \right) \delta \left(
\lambda I-T\right) \,d\lambda =2\pi if\left( T\right) .
\end{gather*}

\begin{example}
Derivation of the orthogonal projection operator onto $\ker \left(
M+I\right) $ of the Hermitian matrix $M$ of the Example \ref{Viola} by the
resolvent technique. \emph{We must compute}%
\begin{equation*}
P_{\lambda _{1}}=\frac{1}{2\pi i}\int_{\left| z+1\right| =1}\mathcal{R}%
\left( z,M\right) \,dz.
\end{equation*}%
\emph{Clearly, we have}%
\begin{equation*}
\mathcal{R}\left( z,M\right) =\frac{1}{z^{2}-z-2}\left[ 
\begin{array}{ccc}
z-1 & 1 & 1 \\ 
1 & z-1 & 1 \\ 
1 & 1 & z-1%
\end{array}%
\right] .
\end{equation*}%
\emph{Using that}%
\begin{equation*}
\int_{\left| z+1\right| =1}\frac{\left\{ 1,z-1\right\} }{\left( z+1\right)
\left( z-2\right) }\,dz=\left\{ -\frac{2\pi i}{3},\frac{4\pi i}{3}\right\}
\end{equation*}%
\emph{we reproduce the result we got earlier.}
\end{example}

\section{The $\protect\delta \left( \protect\lambda I-T\right) $ operator as
a limit}

As $\mu \mapsto \left( \lambda \pm i\epsilon -\mu \right) ^{-1}$ is
continuous, for self-adjoint $T$%
\begin{gather*}
\left( \left( \lambda -i\epsilon \right) I-T\right) ^{-1}-\left( \left(
\lambda +i\epsilon \right) I-T\right) ^{-1} \\
=\int_{-\infty }^{+\infty }\left( \frac{1}{\lambda -i\epsilon -\mu }-\frac{1%
}{\lambda +i\epsilon -\mu }\right) \delta \left( \mu I-T\right) \,d\mu .
\end{gather*}%
If $f\in \mathcal{D}\left( \mathbb{R}\right) $, Proposition \ref{Inter}
yields%
\begin{gather*}
\int_{-\infty }^{+\infty }\frac{f\left( \lambda \right) }{2\pi i}\left\{
\left( \left( \lambda -i\epsilon \right) I-T\right) ^{-1}-\left( \left(
\lambda +i\epsilon \right) I-T\right) ^{-1}\right\} d\lambda \\
=\frac{1}{\pi }\int_{-\infty }^{+\infty }\left\{ \int_{-\infty }^{+\infty
}f\left( \lambda \right) \frac{\epsilon }{\left( \lambda -\mu \right)
^{2}+\epsilon ^{2}}\,\,d\lambda \right\} \delta \left( \mu I-T\right) d\mu .
\end{gather*}%
Since in the sense of distributions%
\begin{equation*}
\frac{1}{2\pi i}\left( \frac{1}{\lambda -i\epsilon -\mu }-\frac{1}{\lambda
+i\epsilon -\mu }\right) \rightarrow \delta \left( \lambda -\mu \right)
\end{equation*}%
as $\epsilon \rightarrow 0^{+}$, we have%
\begin{equation*}
\frac{1}{\pi }\int_{\mathbb{R}}f\left( \lambda \right) \frac{\epsilon }{%
\left( \lambda -\mu \right) ^{2}+\epsilon ^{2}}\,\,d\lambda \rightarrow
\int_{\mathbb{R}}f\left( \lambda \right) \delta \left( \lambda -\mu \right)
d\lambda
\end{equation*}%
as $\epsilon \rightarrow 0^{+}$. Hence, if $g_{n}$ is defined by the
left-hand side\ $\mu $-parametric integral with $\epsilon =1/n$, then $%
g_{n}\rightarrow f$ pointwise on $\mathbb{R}$. So, if $f\in \mathcal{D}%
\left( \mathbb{R}\right) $ and $T$ is bounded (hence with $\sigma \left(
T\right) $ compact), as can be easily checked $\left\{ g_{n}\right\}
_{n=1}^{\infty }$ is a uniformly bounded sequence of continuous functions,
with $\sup_{n\in \mathbb{N}}\left\| g_{n}\right\| _{\infty }\leq \left\|
f\right\| _{\infty }$, that converges pointwise on $\mathbb{R}$\ to $f$.
Thus, by \cite[10.2.8 Corollary]{DS} one has $g_{n}\left( T\right)
\rightarrow f\left( T\right) $ in the strong operator topology, that is%
\begin{gather*}
\frac{1}{\pi }\int_{-\infty }^{+\infty }\left\{ \int_{-\infty }^{+\infty
}f\left( \lambda \right) \frac{\epsilon }{\left( \lambda -\mu \right)
^{2}+\epsilon ^{2}}\,\,d\lambda \right\} \,\delta \left( \mu I-T\right) d\mu
\\
\rightarrow \int_{-\infty }^{+\infty }f\left( \mu \right) \delta \left( \mu
I-T\right) d\mu
\end{gather*}%
as $\epsilon \rightarrow 0^{+}$ in the strong operator topology of $\mathcal{%
L}\left( \mathcal{H}\right) $. Therefore, if $T$ is bounded and $f\in 
\mathcal{D}\left( T\right) $ then%
\begin{equation*}
\int_{-\infty }^{+\infty }\frac{f\left( \lambda \right) }{2\pi i}\left\{
\left( \left( \lambda -i\epsilon \right) I-T\right) ^{-1}-\left( \left(
\lambda +i\epsilon \right) I-T\right) ^{-1}\right\} d\lambda
\end{equation*}%
goes to $f\left( T\right) $ as $\epsilon \rightarrow 0^{+}$. This proves
that for bounded $T$%
\begin{equation*}
\lim_{\epsilon \rightarrow 0^{+}}\frac{1}{2\pi i}\left( \left( \lambda
-i\epsilon \right) I-T\right) ^{-1}-\left( \left( \lambda +i\epsilon \right)
I-T\right) ^{-1}
\end{equation*}%
coincides with $\delta \left( \lambda I-T\right) $ as an $\mathcal{L}\left( 
\mathcal{H}\right) $-valued distribution.

\section{Unitary equivalence of $\protect\delta \left( \protect\lambda %
I-T\right) $}

\begin{theorem}
If $T$ is a self-adjoint operator defined on the whole of $\mathcal{H}$,
there exist a finite measure $\mu $ on the Borel sets of the compact space $%
\sigma \left( T\right) $ and a linear isometry $U$ from $L_{2}\left( \sigma
\left( T\right) ,\mu \right) $ onto $\mathcal{H}$ such that%
\begin{equation*}
U^{-1}\delta \left( \lambda I-T\right) U=\delta \left( \lambda I-Q\right)
\end{equation*}%
where $\left( Q\varphi \right) \left( x\right) =x\,\varphi \left( x\right) $
is the position operator.
\end{theorem}

\proof%
According to \cite{Ha} there exist a finite measure $\mu $ on the Borel sets
of the compact space $\sigma \left( T\right) $ and a linear isometry $U$
from $L_{2}\left( \sigma \left( T\right) ,\mu \right) $ onto $\mathcal{H}$
such that%
\begin{equation*}
\left( U^{-1}TU\right) \varphi =Q\varphi
\end{equation*}%
for every $\varphi \in \mathcal{L}_{2}\left( \sigma \left( T\right) ,\mu
\right) $. So, since $U^{-1}TU$ is a self-adjoint operator on $L_{2}\left(
\sigma \left( T\right) ,\mu \right) $, we have 
\begin{equation*}
U^{-1}\delta \left( \lambda I-T\right) U=\delta \left( \lambda
I-U^{-1}TU\right) =\delta \left( \lambda I-Q\right)
\end{equation*}%
as stated. 
\endproof%

\begin{remark}
For such linear isometry $U$ the equation%
\begin{equation*}
\left( U^{-1}\delta \left( \lambda I-T\right) U\,\varphi \right) \left(
x\right) =\delta \left( \lambda -x\right) \varphi \left( x\right)
\end{equation*}%
holds for every $\varphi \in \mathcal{L}_{2}\left( \sigma \left( T\right)
,\mu \right) $.
\end{remark}

\section{Commutation relations}

Let $S$ and $T$ be two self-adjoint operators defined on the whole of $%
\mathcal{H}$ for which equations $\left[ S,\left[ S,T\right] \right] =\left[
T,\left[ S,T\right] \right] =\mathbf{0}$ hold. In this case%
\begin{equation*}
\left[ -itS,\left[ -itS,-isT\right] \right] =it^{2}s\left[ S,\left[ S,T%
\right] \right] =\mathbf{0}
\end{equation*}%
and the Baker-Campbell-Hausdorff formula yields 
\begin{gather*}
\delta \left( \lambda I-S\right) \delta \left( \mu I-T\right) = \\
\frac{1}{\left( 2\pi \right) ^{2}}\int \int_{\mathbb{R}^{2}}e^{i\left(
t\lambda +s\mu \right) }e^{-itS}e^{-isT}\,dt\,ds= \\
\frac{1}{\left( 2\pi \right) ^{2}}\int \int_{\mathbb{R}^{2}}e^{i\left(
t\lambda +s\mu \right) }e^{\frac{-st}{2}\left[ S,T\right] }e^{-i\left(
tS+sT\right) }\,dt\,ds.
\end{gather*}%
Likewise, since $\left[ T,\left[ T,S\right] \right] =\left[ S,\left[ T,S%
\right] \right] =\mathbf{0}$ one has 
\begin{gather*}
\delta \left( \mu I-T\right) \delta \left( \lambda I-S\right) = \\
\frac{1}{\left( 2\pi \right) ^{2}}\int \int_{\mathbb{R}^{2}}e^{i\left(
t\lambda +s\mu \right) }e^{\frac{-st}{2}\left[ T,S\right] }e^{-i\left(
tT+sS\right) }\,dt\,ds \\
=\frac{1}{\left( 2\pi \right) ^{2}}\int \int_{\mathbb{R}^{2}}e^{i\left(
t\lambda +s\mu \right) }e^{\frac{st}{2}\left[ S,T\right] }e^{-i\left(
tS+sT\right) }\,dt\,ds.
\end{gather*}%
So, using that%
\begin{equation*}
\exp \left( \frac{ist}{2}i\left[ S,T\right] \right) -\left( -\frac{ist}{2}i%
\left[ S,T\right] \right) =2i\sin \left( \frac{ist}{2}\left[ S,T\right]
\right)
\end{equation*}%
we have%
\begin{gather*}
\left[ \delta \left( \lambda I-S\right) ,\delta \left( \mu I-T\right) \right]
= \\
\frac{i}{2\pi ^{2}}\int \int_{\mathbb{R}^{2}}\sin \left( \frac{ist}{2}\left[
S,T\right] \right) e^{i\left( t\lambda +s\mu \right) }e^{-i\left(
tS+sT\right) }\,dt\,ds.
\end{gather*}

For position $Q$ and momentum $P$ of a one-dimensional particle, one has $%
\mathcal{H}=L_{2}\left( \mathbb{R}\right) $ and $\left[ Q,P\right] =i\hbar I$%
. Therefore $\left[ Q,\left[ Q,P\right] \right] =\left[ P,\left[ Q,P\right] %
\right] =\mathbf{0}$ and 
\begin{gather*}
\left[ \delta \left( \lambda I-Q\right) ,\delta \left( \mu I-P\right) \right]
= \\
\frac{i}{2\pi ^{2}}\int \int_{\mathbb{R}^{2}}\sin \left( -\frac{st}{2}\hbar
\right) e^{i\left( t\lambda +s\mu \right) }e^{-i\left( tQ+sP\right)
}\,dt\,ds.
\end{gather*}

According to Theorem \ref{Furie}, if $\left[ \delta \left( \lambda
I-S\right) ,\delta \left( \mu I-T\right) \right] $ acts on $f\left( \lambda
\right) =\lambda $, formally we have%
\begin{gather*}
\int_{-\infty }^{+\infty }\lambda \left[ \delta \left( \lambda I-S\right)
,\delta \left( \mu I-T\right) \right] d\lambda = \\
\frac{i}{2\pi ^{2}}\int \int_{\mathbb{R}^{2}}\sin \left( \frac{ist}{2}\left[
S,T\right] \right) \\
\left\{ \int_{-\infty }^{+\infty }\lambda e^{it\lambda }\,d\lambda \right\}
e^{is\mu }e^{-i\left( tS+sT\right) }\,dt\,ds.
\end{gather*}%
So, using the distributional equality%
\begin{equation}
\int_{-\infty }^{+\infty }\lambda e^{it\lambda }d\lambda =\frac{2\pi }{i}%
\,\delta ^{\prime }\left( t\right)  \label{solt}
\end{equation}%
and integrating by parts, it follows that%
\begin{gather*}
\int_{-\infty }^{+\infty }\lambda \left[ \delta \left( \lambda I-S\right)
,\delta \left( \mu I-T\right) \right] d\lambda \\
=-\frac{i}{2\pi }\left[ S,T\right] \int_{-\infty }^{+\infty }se^{is\left(
\mu I-T\right) }ds.
\end{gather*}%
Observe that a second application of equation (\ref{solt}) and a second
integration by parts yield%
\begin{gather*}
\int_{-\infty }^{+\infty }\lambda \mu \left[ \delta \left( \lambda
I-S\right) ,\delta \left( \mu I-T\right) \right] d\lambda \,d\mu = \\
-\frac{i}{2\pi }\left[ S,T\right] \int_{-\infty }^{+\infty }se^{-isT}\left\{
\int_{-\infty }^{+\infty }\mu e^{is\mu }d\mu \right\} ds= \\
\left[ S,T\right] \int_{-\infty }^{+\infty }\delta \left( s\right) \left\{
1-isT\right\} e^{-isT}ds=\left[ S,T\right]
\end{gather*}%
as expected.

\section{A remark on the Stone formula}

Let $T$ be a self-adjoint operator densely defined on a Hilbert space $%
\mathcal{H}$. If $A$ is a Borel set in $\sigma \left( T\right) $, defining%
\begin{equation}
E\left( A\right) :=\int_{-\infty }^{+\infty }\chi _{A}\left( \lambda \right)
\,dE_{\lambda }  \label{Gladys}
\end{equation}%
where $\chi _{A}$ stands for the characteristic function of $A$ (which is a
bounded Borel function), then $E$ is an $\mathcal{L}\left( \mathcal{H}%
\right) $-valued finitely additive and pointwise countably additive measure
(i.e., countable additivity under the strong operator topology of $\mathcal{L%
}\left( \mathcal{H}\right) $) on the $\sigma $-algebra $\mathcal{A}$ of
Borel subsets of $\sigma \left( T\right) $. So, if the characteristic
function $\chi _{A}$ of $A$ with respect to $\mathbb{R}$ is continuous on $%
\sigma _{p}\left( T\right) $ then 
\begin{equation*}
E\left( A\right) =\int_{-\infty }^{+\infty }\chi _{A}\left( \lambda \right)
\delta \left( \lambda I-T\right) \,d\lambda .
\end{equation*}

For $-\infty <a<b<\infty $ and $\epsilon >0$, we have%
\begin{gather*}
\int_{a}^{b}\left\{ \int_{-\infty }^{+\infty }\left( \frac{1}{\lambda
-i\epsilon -\mu }-\frac{1}{\lambda +i\epsilon -\mu }\right) \delta \left(
\mu I-T\right) \,d\mu \right\} d\lambda \\
=\int_{-\infty }^{+\infty }\left\{ \int_{a}^{b}\frac{2i\epsilon \,d\lambda }{%
\left( \lambda -\mu \right) ^{2}+\epsilon ^{2}}\right\} \,\delta \left( \mu
I-T\right) \,d\mu = \\
2i\int_{\mathbb{R}}\left\{ \arg \tan \left( \!\frac{b-\mu }{\epsilon }%
\!\right) -\arg \tan \left( \!\frac{a-\mu }{\epsilon }\!\right) \right\}
\delta \left( \mu I-T\right) d\mu .
\end{gather*}%
If the limit as $\epsilon \rightarrow 0^{+}$ the bracketed function is equal
to $0$ if $\mu \in \mathbb{R}\setminus \left[ a,b\right] $, equal to $\pi $
if $a<\mu <b$ and equal to $\pi /2$ if $\mu \in \left\{ a,b\right\} $. So,
if $a,b\notin \sigma _{p}\left( T\right) $ so that $\chi _{\left( a,b\right)
}$ and $\chi _{\left[ a,b\right] }$ both belong to $\mathcal{B}_{p}\left( 
\mathbb{R}\right) $, setting%
\begin{equation*}
g_{n}\left( \mu \right) :=\frac{1}{\pi }\int_{a}^{b}\frac{2in^{-1}\,d\lambda 
}{\left( \lambda -\mu \right) ^{2}+n^{-2}}
\end{equation*}%
for each $n\in \mathbb{N}$ and%
\begin{equation*}
f\left( \mu \right) :=\chi _{\left( a,b\right) }\left( \mu \right) +\chi _{
\left[ a,b\right] }\left( \mu \right) ,
\end{equation*}%
then $g_{n}\left( \mu \right) \rightarrow f\left( \mu \right) $ for every $%
\mu \in \mathbb{R}$ and $\sup_{n\in \mathbb{N}}\left\| g_{n}\right\|
_{\infty }\leq 1$ which, according to Proposition \ref{Angela}, implies that 
$g_{n}\left( T\right) x\rightarrow f\left( T\right) x$ for every $x\in
D\left( T\right) $. In other words 
\begin{gather*}
\lim_{\epsilon \rightarrow 0^{+}}\frac{1}{2\pi i}\int_{a}^{b}\left( \left(
\lambda -i\epsilon -T\right) ^{-1}-\left( \lambda +i\epsilon -T\right)
^{-1}\right) \,d\lambda \\
=\frac{1}{2}\int_{-\infty }^{+\infty }\left( \chi _{\left( a,b\right) }+\chi
_{\left[ a,b\right] }\right) \delta \left( \mu I-T\right) \,d\mu ,
\end{gather*}%
holds pointwise on the domain $D\left( T\right) $ of $T$. Hence, by virtue
of (\ref{Gladys}) we get%
\begin{gather*}
\lim_{\epsilon \rightarrow 0^{+}}\frac{1}{2\pi i}\int_{a}^{b}\left( \left(
\lambda -i\epsilon -T\right) ^{-1}-\left( \lambda +i\epsilon -T\right)
^{-1}\right) \,d\lambda = \\
\frac{1}{2}E\left( \left( a,b\right) \right) +\frac{1}{2}E\left( \left[ a,b%
\right] \right) =E\left( a,b\right) +\frac{1}{2}E\left( a\right) +\frac{1}{2}%
E\left( b\right)
\end{gather*}%
which is Stone's formula.

\end{document}